\documentclass[11pt]{article}

\usepackage[T1]{fontenc}
\usepackage{amsmath,amssymb,amstext,dsfont,fancyvrb,float,fontenc,graphicx,subfigure,theorem}
\usepackage{amsfonts,euscript,stmaryrd,enumitem,yhmath,mathrsfs}
 \usepackage{dsfont}
\usepackage{bbm}
\usepackage{color}
\usepackage{hyperref}

\title{\Large\bf From Stochastic Integral {\textit{w.r.t.}} Fractional Brownian Motion to Stochastic Integral {\textit{w.r.t.}} Multifractional Brownian Motion}
\date{}

\cleardoublepage 
\usepackage[strict]{changepage}
\def\abstractname{Abstract -}   
\def\abstract{\begin{adjustwidth}{1cm}{1cm} \par    \footnotesize \noindent {\bf \abstractname}
\def\endabstract{ \end{adjustwidth} \smallskip }}

{\theorembodyfont{\itshape}\newtheorem{theorem}{Theorem}[section]}
{\theorembodyfont{\itshape}\newtheorem{proposition}{Proposition}[section]}
{\theorembodyfont{\itshape}\newtheorem{definition}{Definition}[section]}
{\theorembodyfont{\itshape}\newtheorem{lemma}{Lemma}[section]}
{\theorembodyfont{\itshape}\newtheorem{corollary}{Corollary}[section]}
{\theorembodyfont{\rm}}
{\theorembodyfont{\rm}\newtheorem{remark}{Remark}[section]}
{\theorembodyfont{\rm}}
{\theorembodyfont{\rm }}
\def\proof{{\noindent \bf Proof.~}}

{\theorembodyfont{\itshape}\newtheorem{theodef}{Theorem-Definition}[section]}

\newenvironment{g}[2]{\begin{array}{cl} %
<\hspace{-0.2cm}<\hspace{-0.1cm} #1, \hspace{-0.3cm} &#2 \hspace{-0.1cm} >\hspace{-0.2cm}>}%
{\end{array}}

\newcommand{\h}{\hspace{0.1cm}}

\newcommand{\bD}{\mathbb{D}}

\newcommand{\N}{\mathbf{N}}

\newcommand{\bZ}{\mathbf{Z}}

\newcommand{\sifbm}{\mathbf{B}}

\newcommand{\E}{\mathbf{E}}

\newcommand{\R}{\mathbf{R}}

\newcommand{\cl}{\centerline}

\def\i1{\mathds{1}}  

\def\myqed {{
\parfillskip=0pt 
\widowpenalty=10000 
\displaywidowpenalty=10000 
\leavevmode 
\unskip 
\nobreak 
\hfil 
\penalty50 
\hskip.2em 
\null 
\hfill 
$\square$
\par}} 

\newcommand{\cS}{\mathcal{S}}

\newcommand{\sS}{\mathscr{S}}

\def\cA{{\cal A}}
\def\cB{{\cal B}}

\def\cG{{\mathcal G}}

\def\cS{{\cal S}}

\newcommand{\cF}{\mathcal{F}}

\newcommand{\bit} {\begin{itemize} }
\newcommand{\eit} {\end{itemize} }
\def\ben{\begin{enumerate}}
\def\een{\end{enumerate}}

\def\iti{\item[(i)]}

\def\itii{\item[(ii)]}

\makeatletter

\@addtoreset{equation}{section}
\makeatother

\def\ {\hspace{0.1cm}}

\graphicspath{{./}{figures/}}

\begin{document}
\maketitle
\vspace{-1.25cm}
\cl{{\textit{To appear in Annals of the University of Bucarest}}}

\vspace{0.5cm}

\cl{\textsc{\Large Joachim Lebovits}}
\vspace{0.5cm}

 \begin{abstract}
Because of numerous applications {\textit{e.g.}} in finance and in Internet traffic modelling, stochastic integration w.r.t. fractional Brownian motion (fBm) became a very popular topic in recent years. However, since fBm is not a semi-martingale the It\^{o} integration can not be used for integration w.r.t. fBm and one then needs specific developments.  Multifractional Brownian motion (mBm) is a Gaussian process that generalizes fBm by letting
the local H\"{o}lder exponent vary in time. In addition to the fields mentioned above, it is useful in many and various areas such as geology and biomedicine.
In this work we start from the fact, established in \cite[Thm 2.1.(i)]{fBm_to_mBm_HerbinLebovitsVehel}, that  an mBm may be approximated,
in law, by a sequence of ``tangent" fBms. We used this result to show how one can define a stochastic integral w.r.t. mBm from the  stochastic integral w.r.t. fBm, defined in \cite{Ben1}, in the white noise theory sense.
 \end{abstract}


\begin{MSC}
60G22, 60G15, 60H40, 60H05.
\end{MSC}


\section{Background and motivation}
\label{intro}


\hspace{0.5cm} Fractional Brownian motion (fBm) is a centred Gaussian process that has been used a lot, and still is, these recent years to model many, natural or artificial, phenomena such as physics, geophysics,  financial and Internet traffic modeling, image analysis and synthesis and more.
One of the advantages, that makes the use of fBm so popular to model phenomena, is the fact that fBm is a long range dependence process, that has  the ability to match any prescribed constant local regularity. For any $H$ in $(0,1)$, its covariance function $R_{H}$ reads:
\begin{equation*}
 R_{H}(t,s) :=  \frac{\gamma_{H}}{2} ( {|t|}^{2{H}} + {|s|}^{2{H}} - {|t-s|}^{2{H}}),
\end{equation*}

where $\gamma_{H}$ is a positive constant. The parameter $H$ is called the Hurst index. Besides, when $H=\frac 1 2$, fBm reduces to standard Brownian motion.  Various integral representations of fBm are known, including the harmonizable and moving average ones \cite{TaqSam}, as well as representations by integrals over a finite domain \cite{nualart,dozzi}.

The fact the H\"{o}lder exponent of a fBm remains constant and the same all along its trajectory restricts its application in some situations and more importantly does not seem to be adapted to describe or to model phenomena

 that present in the same time a long range dependance (which requires $H>1/2$) and irregular trajectories (which requires $H<1/2$). Multifractional Brownian motion was introduced to overcome these limitations. The basic idea is to replace the real $H$ by a function $t\mapsto h(t)$, ranging in $(0,1)$.

\smallskip

\hspace{0.5cm}Several definitions of multifractional Brownian motion exist.\h The first ones were proposed in \cite{PL} and in \cite{ABSJDR}. They then have been extended in \cite{StoTaq}. Finally, a more general definition of mBm, that includes all the definitions given in \cite{PL,ABSJDR,StoTaq}, has been introduced  in \cite{fBm_to_mBm_HerbinLebovitsVehel}. In the present work, and for sake of simplicity, we will only deal with  the harmonizable version of mBm, which we will call ``the''  mBm.

More precisely,  a deterministic function $h:\R \rightarrow (0,1)$ being fixed, the mBm of functional parameter $h$, noted $B^h:={(B^h_t)}_{t\in \R}$,  is defined by:

\begin{equation}
\label{ozeijdzoeijdzejdozeij}
B^h_t:= \frac{1}{c_{h(t)}} \ \int_{\R} \frac{e^{itu}-1}{|u|^{h(t)+1/2}} \widetilde{W}(du),
\end{equation}

where $c_x :=  {\big(\frac{2\cos(\pi x) \Gamma(2-2x)}{x(1-2x)} \big)}^{\frac{1}{2}}$, for every $x$ in $(0,1)$,  and where $\widetilde{W}$ denotes a complex-valued Gaussian measure.  Denote $R_{h}$ the covariance function of $B^h$. Thanks to \cite{ACLVL}, we know that:

\vspace{-0.2cm}

\begin{equation*}
R_{h}(t,s)=     \tfrac{ c^2_{h_{t,s}}}{c_{h(t)}c_{h(s)}}  \ \big[\tfrac{1}{2} \big( {|t|}^{2h_{t,s}} + {|s|}^{2h_{t,s}}  - {|t-s|}^{2h_{t,s}} \big)\big],
\end{equation*}

where $h_{t,s} := \frac{h(t)+h(s)}{2}$.

A word on notation: $B^H_.$ or $B^{h(t)}_.$ will always denote an fBm with Hurst index $H$ or $h(t)$, while $B^h_.$ will stand for an mBm.
Note that $B^h_t:= B^{h(t)}_t$, for every real $t$. Besides, and in order to simplify some notations in the sequel, denote $\sifbm(t,H):= \frac{1}{c_{H}} \ \int_{\R} \frac{e^{itu}-1}{|u|^{H+1/2}} \widetilde{W}(du)$, for every $H$ in $(0,1)$.

\subsubsection*{Outline of the paper}
The remaining of this paper is organized as follows. We explain in the first part of Section $2$ how an mBm can be approximated in law by a sequence of ``tangent'' fBms. The second part of Section $2$ is devoted to some heuristics about the way to define a stochastic integral w.r.t.\h mBm as a limit of integrals w.r.t.\h approximating fBms.
Section 3 provides the background on white noise theory and on fractional white noise that will allow us to define, in a rigorous manner, the limiting fractional Wick-It\^{o} integral in Section 4.\h The main result  of Section 4 being Theorem \ref{deftheo}.\h Finally, we compare in Section 5 the limiting fractional Wick-It\^{o} integral to the stochastic integral with respect to mBm that has been defined in \cite{JLJLV1}.

\section{Approximation of mBm and some heuristics}

\hspace{0.5cm}We start this section with the following result, that will be useful in all this paper. Recalling the definition of $\sifbm$ given at the end of the previous section.

\begin{proposition}
\label{frefefe}
 For every  $[a,b]\times [c,d]$  in $\R \times (0,1)$, there exists $\Lambda$ in $\R^*_+$ such that for every  $(t,s,H,H')$ in $[a,b]^2 \times [c,d]^2$,
 \vspace{0.25cm}

 \cl{$\E[{(\sifbm  (t,H) - \sifbm  (s,H'))}^2] \leq \Lambda \ \big( {|t-s|}^{2c} + {|H-H'|}^{2} \big)$,}
\end{proposition}

where $\E[Y]$ denotes the expectation of a real random variable $Y$, which belongs to $L^1(\Omega,\cF,P)$.\\

\proof
Let $(s,t,H,H')$ be fixed in $[a,b]^2\times {[c,d]}^2$. Since, for every $H$ in $[c,d]$, the process  $\{\sifbm  (t,H), \ t \in \R\}$
 is a fractional Brownian motion of Hurst index $H$, we know, using the triangular inequality,  that it is sufficient to show
that $I^{H,H'}_t:=\E[{(\sifbm  (t,H) - \sifbm  (t,H'))}^2] \leq \Lambda \ {|H-H'|}^{2}$. One has:
\vspace{-0.25cm}

$$I^{H,H'}_t  =\int_{\R} \ |\tfrac{e^{it\xi} - 1}{\xi}|^2 \ \left|\tfrac{1}{c_H} \ |\xi|^{1/2-H} - \tfrac{1}{c_{H'}} \ |\xi|^{1/2-H'} \right|^2 \ d\xi.$$

For $\xi$ in $\R^*$, the map $f_{\xi}:[c,d]\rightarrow \R_+$, defined by $f_{\xi}(H):= \tfrac{1}{c_H} \ |\xi|^{1/2-H}$ is $C^1$, since $H\mapsto c_H$ is $C^{\infty}$ on $(0,1)$. Thus there exists a positive real $D$ such that, for all $(\xi,H)$ in $\R^*\times[c,d]$,
\vspace{-0.25cm}

$$ |f'_{\xi}(H)| \leq D \ |\xi|^{1/2-H} \ (1+ |\ln(|\xi|)|) \leq D \ ( |\xi|^{1/2-c} + |\xi|^{1/2-d}) \ (1+ |\ln(|\xi|)|).$$

Thanks to the mean-value theorem, one can write
 \vspace{-0.25cm}

\begin{align*}
I^{H,H'}_t &\leq \ C \ |H-H'|^2 \ \int_{\R} \ \tfrac{|e^{it\xi} - 1|^2}{|\xi|^{2}} \ ( |\xi|^{1/2-c} + |\xi|^{1/2-d} )^2 \ (1+ |\ln(|\xi|)|)^2 \ d\xi  \\
& \leq \ C \ |H-H'|^2 \ \big(   \int_{|\xi| > 1} \ \tfrac{(1+ \ln |\xi| )^2}{|\xi|^{1+2c}} \ d\xi +  \ \ \int_{|\xi| \leq 1} \ |\xi|^{1-2d} \ (1+ |\ln(|\xi|)|)^2 \ d\xi \big)  \notag \\
& \leq C\ |H-H'|^2,
\end{align*}

where $C$  stands for a constant, independent of $t$ and $H$, whose precise value is unimportant and which may change from line to line. \myqed

\begin{remark}
\label{ifedhvreui}

Since  $\sifbm:={(\sifbm(t,H))}_{(t,H)\in \R\times [0,1]}$ is a Gaussian field,  the previous result, as well as  Kolmogorov's criterion, entail that the field $\sifbm$
has a $d$-H\"{o}lder continuous version for any $d$ in $(0,\frac{1}{2} \wedge c)$. In the sequel we will always work with such a version.
\end{remark}

\subsection{Approximation of multifractional Brownian motion}
\label{oidjvofdjvofdivfoi}

\hspace{0.5cm}Since a mBm is generalization of fBm, it is a natural question to wonder if an mBm may be approximated by patching adequately chosen fBms; the sense of this approximation remaining to define. For notational simplicity we take $[a,b]:=[0,1]$ in the sequel.
Heuristically,  we divide $[0,1)$ into ``small'' intervals $[t_i, t_{i+1})$, and replace on each of these $B^h$ by the fBm $B^{H_i}$ where $H_i = h(t_i)$.
It seems reasonable to expect that the resulting process $\sum_i B^{H_i}_t \mathbbm{1}_{[t_i, t_{i+1})}(t)$ will converge, in law, to $B^h$ when the sizes of the intervals $[t_i, t_{i+1})$ go to 0. The notations remains the same as previously. Let $h:[0,1]\rightarrow (0,1)$ be a conitnuous deterministic function and denote $B^h$ the fixed mBm of functional parameter $h$.
 Let us explain how  this mBm can be approximated on $[0,1]$ by patching together fractional Brownian motions defined on a sequence of partitions of $[0,1]$.

\vspace{0.25cm}

In that view,  we choose an increasing sequence $(q_n)_{n\in \N}$ of integers such that $q_0:=1$ and $2^n \leq q_n \leq2^{2^n}$ for all $n$ in $\N^*$.
For any $n$ in $\N$, define $x^{(n)}:=\{ x^{(n)}_{k}; k\in [\hspace{-0.05cm}[ 0,q_n ]\hspace{-0.05cm}] \}$ where $ x^{(n)}_{k}:= \frac{k}{q_n}$ for $k$ in
$[\hspace{-0.05cm}[ 0, q_n ]\hspace{-0.05cm}]$
(for integers $p$ and $q$ with $p<q$, $[\hspace{-0.05cm}[ p, q ]\hspace{-0.05cm}]$ denotes the set $\{p;p+1;\cdots;q\}$). Define, for  $n$ in $\N$, the partition $\cA_n:=\{ [x^{(n)}_{k}, x^{(n)}_{k+1} ); k\in [\hspace{-0.05cm}[ 0, q_n -1 ]\hspace{-0.05cm}]\} \cup \{x^{(n)}_{q_n}\}$. It is clear that $\cA:=(\cA_n)_{n\in\N}$ is a decreasing nested sequence of subdivisions of $[a,b]$ (\textit{i.e.} $\cA_{n+1} \subset  \cA_n$, for every $n$ in $\N$).

For $t$ in $[0,1]$ and $n$ in $\N$ there exists a unique integer $p$ in $[\hspace{-0.05cm}[ 0, q_n-1 ]\hspace{-0.05cm}]$ such that $x^{(n)}_p \leq t <  x^{(n)}_{p+1}$. We will note $x^{(n)}_t$ the real $x^{(n)}_p$ in the sequel.  It is clear that the sequence ${(x^{(n)}_t)}_{n\in\N}$ is increasing and converges to $t$ as $n$ tends to $+\infty$. Besides, define for $n$ in $\N$, the function $h_n:[0,1] \rightarrow (0,1)$ by setting $h_n(1) = h(1)$ and, for any $t$ in $[0,1)$, $h_n(t):= h(x^{(n)}_t)$.
The sequence of step functions ${(h_n)}_{n\in \N}$ converges  pointwise to $h$ on $[0,1]$. Define, for $t$ in $[0,1]$ and $n$ in $\N$, the process

\vspace{-0.5cm}

\begin{equation}
\label{eq:Xn}
B^{h_n}_t := \sifbm  (t,h_n(t)) =  \sum_{k=0}^{q_n - 1} \mathbbm{1}_{[x^{(n)}_k, x^{(n)}_{k+1})}(t)\ \sifbm  (t,h(x^{(n)}_k)) + \mathbbm{1}_{\{1\}}(t)\  \sifbm  (1,h(1)) .
\end{equation}


Note that, despite the notation, the process $B^{h_n}$ is not an mBm, as $h_n$ is not continuous.

The following theorem, the proof of which can be found in \cite[Thm 2.1.(i)]{fBm_to_mBm_HerbinLebovitsVehel}, shows that mBm appears naturally as a limit object of sums of fBms.
%

{\begin{theorem}[Approximation theorem]
\label{th:approxMBMweak}

Let $\cA$ be a sequence of partitions of $[0,1]$ as defined above, and consider the sequence of processes defined in \eqref{eq:Xn}.

 If $h$ is $\beta$-H\"{o}lder continuous for some positive real $\beta$, then the sequence of processes $(B^{h_n})_{n\in\N}$  converges, in law, to the process $B^h$. In oher words,

 $$\{B^{h_n}_t; \ {t\in [a,b]} \}\xrightarrow[n \to +\infty]{\text{law}}   \{B^{h}_t; \ {t\in [a,b]} \}.$$
\end{theorem}
}

\subsection{Some Heuristics about Stochastic integrals w.r.t. mBm as limits of integrals w.r.t. fBm}
\label{oerifjeorifjeroifjeroifjeo}
\hspace{0.5cm}In the remaining of this paper, we consider a $C^1$ deterministic function $h:\R \rightarrow (0,1)$.
Moreover $\sifbm$ still denotes the Gaussian field defined at the end of Section 1 and $B^{h}$ the mBm defined by \eqref{ozeijdzoeijdzejdozeij}. The mBm $B^h$  then verifies $B^h_t\stackrel{\text{a.s.}}{:=} \sifbm(t,h(t))$, for every $t$.

The result of Theorem \ref{th:approxMBMweak} suggests that one may define stochastic integrals with respect to mBm as limits of integrals with respect to approximating fBms.  Let first notice the following fact, the proof of which can be found in \cite[Appendix B]{fBm_to_mBm_HerbinLebovitsVehel}, that will be useful to define the integral with respect to mBm.

\begin{proposition}
 \label{idazuhdiz}
  For every real $t$, the map $H\mapsto \sifbm(t,H)$, from $(0,1)$ to $L^2(\Omega)$, is $C^1$. Moreover, it fulfills the following H\"{o}lder condition.
  For all  $[a,b] \times [c,d] \subset \R \times (0,1)$ there exists $\Delta$ in $\R^*_+$ such that:

$$\E\big[{(\tfrac{\partial \sifbm}{\partial H}(t,H) - \tfrac{\partial \sifbm}{\partial H}(s,H'))}^2\big] \leq \Delta \ ({|t-s|}^{2} + {|H-H'|}^{2}),$$

where $\tfrac{\partial \sifbm}{\partial H}(t,H)$ denotes, for every real $t$, the derivative with respect to $H$ of the map $H\mapsto \sifbm(t,H)$.
 \end{proposition}

Let us explain now, in a heuristic way, how to define an integral with respect to mBm using approximating fBms. Write the ``differential'' of
$\sifbm(t,H)$:

\begin{equation*}
d\sifbm(t,H) = \frac{\partial \sifbm}{\partial t }(t,H) \ dt + \frac{\partial \sifbm}{\partial H }(t,H) \ dH.
\end{equation*}

Of course, this is only formal as $t \mapsto\sifbm(t,H)$ is not differentiable in the $L^2$-sense nor almost surely with respect to $t$. It is, however, in the sense of Hida distributions (this will be made precise in Section \ref{whitenoise}, see in particular \eqref{fiizeufhirsozeedeijdefrt2}). With a differentiable function $h$ in place of $H$, this (again formally) yields

\begin{equation}
\label{frejo}
d\sifbm(t,h(t))  =  \frac{\partial \sifbm}{\partial t }(t,h(t)) \ dt +  h'(t) \ \frac{\partial \sifbm}{\partial H }(t,h(t))  \ dt.
\end{equation}

Of course the first term of the right hand side of \eqref{frejo} has no meaning \textit{a priori} since mBm is not differentiable with respect to $t$. However, continuing with our heuristic reasoning, we then approximate $\frac{\partial \sifbm}{\partial t }(t,h(t))$ by $\lim\limits_{n\to +\infty}  \sum^{q_n-1}_{k=0} \mathbbm{1}_{[x^{(n)}_k, x^{(n)}_{k+1})}(t)\ \frac{\partial \sifbm}{\partial t }(t,h_n(t))$. This formally yields:
\vspace{-0.25cm}

\begin{equation}
\label{frejoeedeezzedezdzeezd}
d\sifbm(t,h(t))    \approx  {\textstyle{\lim\limits_{n\to +\infty}  \sum^{q_n-1}_{k=0}} \mathbbm{1}_{[x^{(n)}_k, x^{(n)}_{k+1})}(t)} \ \tfrac{\partial \sifbm}{\partial t }(t,h_n(t)) \ dt +  h'(t) \ {\tfrac{\partial \sifbm}{\partial H }(t,h(t)) }  \ dt.
\end{equation}

For the sake of notational simplicity, we will consider integrals over the interval $[0,1]$.
Let us note $\int_{[0,1]} Y_t \ d^{\diamond}B^{H}_t$ the integral of a process Y with respect to a fBm of Hurst index $H$, in the white noise theory sense (that will be fully detailed in the next section), assuming it exists.

From \eqref{frejoeedeezzedezdzeezd}, and assuming we may exchange integrals and limits, we would thus like to define, for suitable processes $Y$, the integral with respect to the mBm $B^h$, noted $\int^1_0\ Y_t \ d^{\diamond}B^h_t$,  by setting,

\vspace{-0.5cm}

\begin{align}
\label{frejoeededd}
{\textstyle \int^1_0\ Y_t \ d^{\diamond}B^h_t=  \lim\limits_{n\to +\infty} \sum^{q_n-1}_{k=0} \  \int^{x^{(n)}_{k+1}}_{x^{(n)}_k} \ Y_t \  d^{\diamond}B^{h(x^{(n)}_k)}_t   + \int^1_0 \ Y_t \   h'(t) \ {\frac{\partial \sifbm}{\partial H }(t,h(t)) }  \ dt},
\end{align}

where the first term of the right-hand side of \eqref{frejoeededd} is a limit, in a sense to be made precise
 of a sum of integrals with respect to fBms and the second term is a weak integral (see Section \ref{whitenoise}). The following notation will be useful:
\medskip

{\bf Notation (integral with respect to lumped fBms)}
{\it
\label{isudhcidsuhcvuhsihiudchs}
Let  $Y:={(Y_t)}_{t\in[0,1]}$ be a real-valued process on $[0,1]$ which is integrable with respect to all fBms of index $H$ in $h([0,1])$ in the white noise theory sense. We denote the integral with respect to lumped fBms  by:
\begin{equation}
\label{eq:Xndede}
\int^1_0 \ Y_t \ d^{\diamond}B^{h_n}_t :=   \sum_{k=0}^{q_n - 1} \int^1_0  \  \mathbbm{1}_{[x^{(n)}_k, x^{(n)}_{k+1})}(t)\ \ Y_t \ d^{\diamond}B^{h(x^{(n)}_k)}_t, \hspace{0.35cm}  n \in \N,
\end{equation}
}

where  $(q_n)_{n\in\N}$ and $x^{(n)}:=\{ x^{(n)}_{k}; k\in [\hspace{-0.05cm}[ 0,q_n ]\hspace{-0.05cm}] \}$ have been defined in Section \ref{oidjvofdjvofdivfoi}. With this notation, our tentative definition of an integral w.r.t. to mBm \eqref{frejoeededd} reads:
\begin{equation*}
\int^1_0\ Y_t \ d^{\diamond}B^h_t  =  \lim\limits_{n\to +\infty} \  \int^1_0 \ Y_t \ d^{\diamond}B^{h_n}_t   + \int^1_0 \ Y_t \   h'(t) \ {\tfrac{\partial \sifbm}{\partial H }(t,h(t)) }  \ dt,
\end{equation*}

The drawback with the previous definition is that, when $\int^1_0 \ Y_t \   h'(t) \ {\frac{\partial \sifbm}{\partial H }(t,h(t)) }  \ dt
$ will belong to $L^2(\Omega)$, it will not be centred, {\textit{a priori}}. For this reason and because one can see the Wick product as a centered product, we would rather choose, as a definition of $\int^1_0\ Y_t \ d^{\diamond}B^h_t $, the following one:

\vspace{-0.5cm}

\begin{align}
\label{tentativedef}
\int^1_0\ Y_t \ d^{\diamond}B^h_t  &:=  \lim\limits_{n\to +\infty} \  \int^1_0 \ Y_t \ d^{\diamond}B^{h_n}_t   + \int^1_0 \ h'(t) \  Y_t \  \diamond {\tfrac{\partial \sifbm}{\partial H }(t,h(t)) }  \ dt,
\end{align}

where $\diamond$ denotes the Wick product (that will be rigorously defined in the next section).

In order to make the previous statement rigorous, we need to give a precise meaning to the right hand side of \eqref{tentativedef}. In particular, giving a precise meaning to $\int^1_0 \ Y_t \ d^{\diamond}B^{H}_t$ and thus to $\int^1_0 \ Y_t \ d^{\diamond}B^{h_n}_t$ is crucial. This is the aim of the next section.

\section{Backgrounds on white noise theory and Fractional Wick-It\^{o} Integral in the Bochner sense}
\label{whitenoise}

This section is divided into two parts.
In the first one we briefly recall some basic facts about white noise theory and the Bochner integral.\h
In the second part we particularize the definition of the Fractional Wick-It\^{o} Integral, defined in  \cite{Ben1,Ben2,bosw,ell}, into the framewok of Bochner integral.

\subsection{Recalls on white noise theory and the Bochner integral}\label{WNB}
\subsubsection{White noise theory}
\label{suhidshicudiuchdsiuhcsdihdshiu}
%

The following subsection being on purpose extremely short. The reader who is no familiar with white noise theory should refer to \cite{Kuo2} and references therein.

Define the measurable space $(\Omega,\cF)$ by setting  $\Omega:= {\sS}^{'}(\R)$ and $\cF := \cB({\sS}^{'}(\R))$, where $\cB$ denotes the $\sigma$-algebra of Borel sets. Denotes $\mu$ the unique probability measure on $(\Omega,\cF)$ such that, for every $f$ in $L^2(\R)$, the map $<.,f>:\Omega \rightarrow \R$ defined by $<.,f>(\omega) = <\omega,f>$ (where $< , >$ continuously in $L^2(\R)$ extends the action of tempered distributions on Schwartz functions)
is a centred Gaussian random variable with variance equal to ${\|f\|}^2_{L^2(\R)}$ under $\mu$.

We also denote $(L^2)$ the space $L^2(\Omega,\cG,\mu)$ where $\cG$ is the $\sigma$-field generated by ${(<.,f>)}_{f \in L^2(\R)}$, and for every $n$ in $\N$,  define  the $n-$th Hermite function by $e_n(x):=   {(-1)}^n \ {\pi}^{-1/4} {(2^n n!)}^{-1/2} e^{x^2/2} \frac{d^n}{dx^n}( e^{-x^2})$. Denote $A$ the operator defined on $\sS(\R)$ by $A:= -\frac{d^2}{dx^2} + x^2 +1$ and $\Gamma(A)$ the second quantization operator of $A$ (see \cite[Section 4.2]{Kuo2}).

Denote, for $\varphi$ in $(L^2)$, ${\|\varphi\|}^2_0:={\|\varphi\|}^2_{(L^2)}$  and, for $n$ in $\N$, let $\bD\text{om}({\Gamma(A)}^n)$ be the domain of the  $n-$th iteration of $\Gamma(A)$. Define the  family  of norms  ${({\|\ \|}_p)}_{p \in \bZ}$ by:

\vspace{-0.5cm}

\begin{equation*}
  {\|\Phi\|}_p :=  {\|\Gamma(A)^p\Phi\|}_{0} =  {\|\Gamma(A)^p\Phi\|}_{(L^2)},  \h \forall p \in \bZ,\hspace{0.25cm} \forall\Phi \in (L^2)\cap \bD{\text{om}}({\Gamma(A)}^p).
\end{equation*}

For $p$ in $\N$, define $({\cS}_{p}):=\{\Phi \in (L^2): \  \Gamma(A)^p\Phi \ \text{exists and belongs to}  \ (L^2) \}$ and define
$({\cS}_{-p})$ as the completion of the space  $(L^2)$ with respect to the norm ${{\|\ \|}_{-p}}$.
As in \cite{Kuo2}, we let $(\cS)$ denote the projective limit of the sequence ${( (\cS_{p}))}_{p \in \N}$ and ${(\cS)}^*$ the inductive limit of the sequence  ${(({\cS_{-p}}))}_{p \in \N}$. This means that we have the equalities $(\cS) = \underset{p \in \N }{\cap}{({\cS}_{p})}$ (resp.  ${(\cS)}^*= \underset{p \in \N }{\cup}{({\cS_{-p}})}$) and  that convergence in $(\cS)$ (resp. in ${(\cS)}^*$) means convergence in $({\cS}_{p})$ for every $p$ in $\N$   (resp. convergence in  $({\cS}_{-p})$ for some $p$ in $\N$ ).

The space  $(\cS)$ is called the space of stochastic test functions and ${(\cS)}^*$  the space of Hida distributions. Since ${(\cS)}^*$ is the dual space of $(\cS)$. We will note  $<\hspace{-0.2cm}<\hspace{-0.1cm} \ \ , \  \hspace{-0.1cm}>\hspace{-0.2cm}>$  the duality bracket between ${(\cS)}^*$ and $(\cS)$. If $\phi$ and $\Phi$ both belong to $(L^2)$ then we have the equality $\begin{g}{\Phi}{\varphi}\end{g} = {<\Phi,\varphi>}_{(L^2)} = \E[\Phi \ \varphi]$.

A function $\Phi:\R \rightarrow {(\cS)}^{*}$ is called a stochastic distribution process, or an ${(\cS)}^{*}\hspace{-0.1cm}-$process, or a Hida
process. A Hida process $\Phi$ is said to be differentiable at $t_0 \in \R$ if $\lim\limits_{r \to 0} \ r^{-1} ( \Phi(t_0+r)-\Phi(t_0))$ exists in ${(\cS)}^{*}$.

The $S$-transform of an element $\Phi$ of $(\cS^*)$, noted $S(\Phi)$, is defined as the function from  $\sS(\R)$ to $\R$ given,  for every $\eta$ in $\sS(\R)$, by $S(\Phi)(\eta):= {\begin{g}{\Phi}{e^{<.,f> - \frac{1}{2}{|f|}^2_0}\hspace{0.075cm}}  \end{g} }$, where
 ${({|\ \ |}_p)}_{p\in\bZ}$ is  the family  norms defined by ${|f|}^2_p:=   \sum^{+\infty}_{k=0} {(2k + 2)}^{2p} \ {<f,e_k>}^2_{L^2(\R)}$, for all $(p,f)$ in $\bZ \times L^2(\R)$.

Finally for every $(\Phi,\Psi)\in {(\cS)}^{*}\times {(\cS)}^{*}$, there exists a unique element of ${(\cS)}^{*}$, called the Wick product of  $\Phi$ and $\Psi$  and  noted $\Phi \diamond \Psi$,  such that $S(\Phi\diamond \Psi)(\eta) = S(\Phi)(\eta) \ S(\Psi)(\eta)$ for every $\eta$ in $\sS(\R)$.

\subsubsection{Fractional and multifractional White noise}

We now introduce two operators, denoted  $M_H$ and $\tfrac{\partial M_H}{\partial H }$, that will prove useful for the definition of the integral with respect to fBm and mBm.

\subsubsection*{Operators $M_H$ and $\tfrac{\partial M_H}{\partial H }$}
\label{ksksksks}

Let $H$ be a fixed real in $(0,1)$. Following \cite{ell} and references therein, define the operator $M_H$, specified in the Fourier domain, by $\widehat{M_H(u)}(y) := \tfrac{\sqrt{2\pi}}{c_H} \hspace{0.1cm}{|y|}^{1/2-H} \widehat{u}(y)$ for every $y$ in $\R^*$.
 This operator is well defined on the homogeneous Sobolev space of order $1/2-H$, denoted $L^2_H(\R)$ and defined by $L^2_H(\R):=  \{u \in
 \mathcal{\sS'(\R)} \hspace{0.1cm}:\hspace{0.1cm} \widehat{u}  \in L^{1}_{loc}(\R) \hspace{0.1cm}\text{ and }\hspace{0.1cm}
 {\|u\|}_{H} <  +\infty  \}$, where the norm ${\|\ \|}_{H}$ derives from the inner product ${< , >}_{H}$ defined on $L^2_H(\R)$ by ${<u,v>}_{H} := \frac{1}
 {c^2_H} \int_{\R} {|\xi|}^{1-2H} {\widehat{u \hspace{0.1cm}}(\xi) } \ \overline{{\widehat{v \hspace{0.1cm}}(\xi) }} \  d{\xi}$ and where $c_{H}$ was given in
\eqref{ozeijdzoeijdzejdozeij}.

The definition of the operator  $\tfrac{\partial M_H}{\partial H }$ is quite similar. More precisely, define for every $H$ in $(0,1)$, the space $\Gamma_H(\R) := \{u \in \mathcal{\sS'(\R)} \hspace{0.1cm} : \hspace{0.1cm} \widehat{u} \in L^{1}_{loc}(\R) \hspace{0.1cm}\text{  and}\hspace{0.1cm}  {\|u\|}_{\delta_H(\R)} <  +\infty  \}$, where the norm ${\|\ \|}_{\delta_H(\R)}$ derives from the inner product on $\Gamma_H(\R)$ defined by ${<u,v>}_{\delta_H} := \frac{1}{c^2_H} \int_{\R} {(\beta_H + \ln|\xi|)}^2 \ {|\xi|}^{1-2H}   \ {\widehat{u \hspace{0.1cm}}(\xi) } \ \overline{{\widehat{v \hspace{0.1cm}}(\xi) }}  \ d{\xi}$.
Following \cite{JLJLV1}, define the operator $\tfrac{\partial M_H}{\partial H }$ from  $(\Gamma_H(\R), {<,>}_{\delta_H(\R)})$ to  $(L^2(\R), {<,>}_{L^2(\R)})$, in the Fourier domain, by: $\widehat{\tfrac{\partial M_H}{\partial H }(u)}(y) :=  -(\beta_H + \ln|y|) \ \tfrac{\sqrt{2\pi}}{c_H}  \ {|y|}^{1/2 - H}  \ \widehat{u}(y) $, for every  $y$ in $\R^*$.
The reader interested in the properties of $M_H$ and $\tfrac{\partial M_H}{\partial H }$ may refer to \cite[Sections $2.2$ and $4.2$]{JLJLV1}.

\subsubsection*{Fractional and multifractional White noise}

Recall the following result (\cite[(5.10)]{JLJLV1}): Almost surely, for every $t$,

\begin{equation}
\label{pzokpdskrefopsodfkpodkf}
 B^h_t   = { \sum^{+ \infty}_{k = 0} \big( {\textstyle \int^t_0}  \ M_{h(t)}({e}_k)(s) \ ds \big) <.,e_k>.}
\end{equation}

We now define the derivative in the sense of $(\cS)^*$ of mBm. Define the $(\cS^*)$-valued process $W^h:={(W^h_t)}_{t\in[0,1]}$  by

\vspace{-0.25cm}
\begin{equation}
\label{fiizeufhirsozeedeijdefrt2}
W^h_t := {  \sum^{+ \infty}_{k = 0} \big[\tfrac{d}{dt} \big({\textstyle \int^t_0}  \ M_{h(t)}({e}_k)(s)\ ds \big) \big]  \ <.,e_k>.}
\end{equation}

\begin{theodef}\cite[Theorem-definition 5.1]{JLJLV1}
\label{alalalalalal2}
The process $W^{h}$ defined by $(\ref{fiizeufhirsozeedeijdefrt2})$ is an ${(\cS)}^*$-process which verifies, in ${(\cS)}^*$, the following equality:

\vspace{-0.5cm}

\begin{equation}
\label{firsozoijkqsjdefrt2}
 W^h_t =    \sum^{+ \infty}_{k = 0}   M_{h(t)}({e}_k)(t)  <.,e_k> + \ h'(t)\ \sum^{+ \infty}_{k = 0} \big( {\textstyle \int^t_0 }\
\tfrac{\partial M_{H}}{\partial H}({e}_k) (s) \big|_{H=h(t)}  ds \big) <.,e_k>  .
\end{equation}

Moreover the process $B^{h}$ is ${(\cS)}^*$-differentiable on $[0,1]$ and verifies $\frac{dB^{h}}{dt}(t) = W^h_t$ in ${(\cS)}^*$.
\end{theodef}

When the function $h$ is constant and identically equal to $H$, we will write  $W^H:={(W^H_t)}_{t\in [0,1]}$ and call the ${(\cS)}^*$-process $W^{H}$ a fractional white noise. Note that \eqref{firsozoijkqsjdefrt2} may be written as

\begin{equation}
\label{firsozoijkqdezdzedezddesjdefrt}
 W^h_t =   W^{h(t)}_t + h'(t) \ \tfrac{\partial \sifbm}{\partial H}(t,h(t)),
\end{equation}

where $ W^{h(t)}_t$ is nothing but  $W^{H}_t|_{H = h(t)}$ and where the equality holds in ${(\cS)}^*$.

\subsubsection{Bochner integral}
\label{appendiceB2}

Since the objects we are dealing with are no longer random variables in general, the Riemann or Lebesgue integrals are not relevant here. However, taking advantage of the fact that we are working with vector linear spaces, we may use Pettis or Bochner integrals.
In the frame of the Wick-It\^{o} integral, and in view of the result that will provided by Lemma \ref{szfojiojfdsfoifoj} below, the use of Bochner integral appears to be relevant. Indeed, Lemma \ref{szfojiojfdsfoifoj}  gives an easy criterion for integrability, w.r.t. fBm, of any ${(\cS)}^*$-valued process $Y$.
Thus we give here a brief statement on Bochner integral. However, and in order not to weigh down this statement we will only give the necessary tools to proceed (see \cite[p.$247$]{Kuo2} for more details about Bochner integral).

\begin{definition}[Bochner integral \cite{Kuo2}, p.$247$]
\label{bb2}
Let $I$ be a Borel subset of $[0,1]$ and $\Phi:={(\Phi_t)}_{t\in I}$ be an ${(\cS)}^*$-valued process verifying:
\bit
\iti the process $\Phi$ is weakly measurable on $I$ {\textit{i.e.}} the map $t \mapsto {<\hspace{-0.2cm}<\hspace{-0.1cm} \ \Phi_t ,\varphi \ \hspace{-0.1cm}>\hspace{-0.2cm}>}$ is measurable on $I$, for every $\varphi$ in ${(\cS)}$.
\itii there exists $p \in \N$ such that $\Phi_t \in ({\cS}_{-p})$ for almost every $t \in I$ and $t \mapsto {\|\Phi_t\|}_{-p}$ belongs to $L^1(I)$.
\label{foidfj}
\eit

Then there exists an unique element in ${(\cS)}^{*}$, noted $\int_{I} \Phi_{u} \ du $, called the Bochner integral of $\Phi$ on $I$
such that,  for all  $\varphi$ in $(\cS)$,

\vspace{-0.25cm}

$${<\hspace{-0.2cm}<\hspace{-0.1cm} \ \int_{I} \Phi_u \ du,\varphi \ \hspace{-0.1cm}>\hspace{-0.2cm}>} =
\int_{I}
 <\hspace{-0.2cm}<\hspace{-0.1cm} \ \Phi_{u},\varphi\  \hspace{-0.1cm}>\hspace{-0.2cm}> \ du.$$

In this latter case one says that $\Phi$ is Bochner-integrable on $I$ with index $p$.

\end{definition}

\begin{proposition}
\label{oiefhierhiuerhieru}
If $\Phi\hspace{-0.1cm} :I \hspace{-0.15cm} \rightarrow \hspace{-0.1cm} {(\cS)}^{*}$ is Bochner-integrable on $I$ with
index $p$ then ${\| \int_{I} \Phi_t  \ dt \|}_{-p} \leq \int_{I}  {\|\Phi_t\|}_{-p} \ dt.$
\end{proposition}

\begin{theorem}[\cite{Kuo2}, Theorem 13.5]
\label{ccdezd}
Let $\Phi:={(\Phi_t)}_{t\in[0,1]}$ be an ${(\cS)}^{*}$-valued process such that:
\bit
\iti $t\mapsto S(\Phi_t)(\eta)$ is measurable for every $\eta$ in $\sS(\R)$.
\itii There exist $p$ in $\N$, $b$ in $\R^+$ and a function $L$ in $L^1([0,1], dt)$ such that, for \textit{a.e.} $t$ in $[0,1]$, $|S(\Phi_t)(\eta)| \leq L(t) \ e^{b {| \eta |}^2_{p} }$, for every $\eta$ in $\sS(\R)$.
\eit
Then $\Phi$ is Bochner integrable on $[0,1]$  and $\int^1_0 \Phi(s) \ ds \in ({\cS}_{-q})$ for every $q > p$ such that $2 b e^2 D(q-p) < 1$ where $e$ denotes the base of the natural logarithm and where $D(r):= \frac{1}{2^{2r}} \ \sum^{+\infty}_{n=1} \frac{1}{n^{2r}}$ for $r$ in $(1/2,+\infty)$.
\end{theorem}

\noindent
\subsection{Wick-It\^{o} integral with respect to fBm in the Bochner sense}
\label{dsodpqsldqpslpdsqplqdsplqdslsdlqspdlqspdlsqplpdfogdfjfygdufysguy}

The fractional Wick-It\^{o} integral with respect to fBm (or integral w.r.t. fBm in the white noise sense) was introduced in \cite{ell} and extended in \cite{Ben1} using the Pettis integral. As we will see in Lemma \ref{szfojiojfdsfoifoj} below, the Bochner integrability of an ${(\cS)}^*$-valued process $Y$ is a simple condition that ensures the Wick-It\^{o} integrability of $Y$ with respect to fBm (see Definition \ref{oezifhherioiheroiuh2}
 below) of any Hurst index $H$ in $(0,1)$.
For this reason, we now particularize  the fractional Wick-It\^{o} integral with respect to fBm (or Wick-It\^{o} integral w.r.t. fBm)  of \cite{ell} and  \cite{Ben1} in the framework of the Bochner integral. 

\begin{definition}[Wick-It\^{o} integral w.r.t fBm in the Bochner sense]
\label{oezifhherioiheroiuh2}
Let $H \in (0,1)$, $I$ be a Borel subset of $[0,1]$, $B^{H}:={(B^H_t)}_{t\in I}$ be a fractional Brownian motion of Hurst index $H$, and $Y:={(Y_t)}_{t\in I}$ be an ${(\cS)}^*$-valued process verifying:
\bit
\iti there exists $p \in \N$ such that $Y_t \in ({\cS}_{-p})$ for almost every $t \in I$,
\label{foidfj}
\itii the process $t \mapsto Y_t\diamond W^{H}_t$ is Bochner integrable on $I$.
\eit

Then, $Y$ is said to be Bochner-integrable with respect to fBm on $I$ and its integral is defined by:

\begin{equation}
\label{eigfrretth2}
\int_{I} Y_s\ d^{\diamond}B^{H}_s :=  \int_{I} Y_s \diamond W^{H}_s ds.
\end{equation}
\end{definition}

\begin{remark}
\label{pferfoeroferjfoije}
In order to keep the name given in \cite{ell}, we also call this integral fractional Wick-It\^{o} integral.
\end{remark}

The following lemma ensures us that every Bochner integrable process is integrable on $[0,1]$ w.r.t.\h fBm of any  Hurst index $H$ in $(0,1)$. For sake of notational symplicity one denotes, for every integer $p_0$, $q(p_0)$ the integer equal to $\max\{p_0 +1;3\}$ if $p_0 \geq 1$ and equal to $2$ if $p_0=0$.

\begin{lemma}
\label{szfojiojfdsfoifoj}
Let $Y:={(Y_t)}_{t\in[0,1]}$ be an ${(\cS)}^*$-valued process, Bochner integrable of index $p_0 \in \N$. Then $Y$ is integrable on $[0,1]$, with respect to fBm of any Hurst index $H$, in the Bochner sense. Moreover, for any $H$ in $(0,1)$,  $\int_{[0,1]} Y_s\ d^{\diamond}B^{H}_s$ belongs to $({\cS}_{-q(p_0)})$.
\end{lemma}

\proof
Fix $H \in (0,1)$ and an integer $p_0 \geq 2$.
 The map $t\mapsto Y_t \diamond W^{H}_t$ is weakly measurable since $t\mapsto S(Y_t \diamond W^{H}_t)(\eta)$ is measurable for all $\eta$ in $\sS
 (\R)$. Using \cite[Remark $2$ p.$92$]{Kuo2}, one obtains that, for almost all $t$ in $[0,1]$,  ${\|Y_t \diamond W^{H}_t\|}_{-q(p_0)} \leq {\| Y_t \|}_{-p_0}   \
{\| W^{H}_t \|}_{-p_0}$. Hence $Y_t \diamond W^{H}_t$ belongs to $({\cS}_{-q(p_0)})$. Since the map $t\mapsto {\|W^{H}_t\|}_{-r}$ is continuous for all integer $r\geq 2$ (see  \cite[Proposition 5.9]{JLJLV1}), one also gets:

\begin{equation*}
 {\textstyle \int^1_0 } {\|Y_t \diamond W^{H}_t\|}_{-q(p_0)} \    dt     \leq  (\underset{t \in [0,1] }{\sup}{{\|W^{H}_t\|}_{-p_0}}) \  {\textstyle \int^1_0}  {\|Y_t\|}_{-p_0} \ dt <+\infty.
\end{equation*}

This shows that $t \mapsto Y_t \diamond W^{H}_t$ is  Bochner-integrable of index $q(p_0)$.

Let us now assume that $p_0$ belongs to $\{0,1\}$. It is sufficient to check that Theorem  \ref{ccdezd} applies. Condition $(i)$ is obviously fulfilled.
Moreover, using \cite[p.$79$]{Kuo2}, we obtain that, for every $(t,\eta)$ in $[0,1]\times\sS(\R)$,

\begin{equation*}
|S(Y_t \diamond W^{H}_t)(\eta)| \leq \|Y_t\|_{-p_0} \ \underset{t \in [0,1] }{\sup}{{\|W^{H}_t\|}_{-2}} \ e^{ \frac{1}{2}{| \eta |}^2_{2} } =: L(t) \ e^{ \frac{1}{2}{| \eta |}^2_{2} }.
\end{equation*}

Since $Y$ is Bochner integrable of index $p_0$, it is clear that $L$ belongs to $L^1([0,1],dt)$. Moreover, $e^2 D(r-p_0) < 1$, for every  $r \geq p_0 +2$. Theorem  \ref{ccdezd} then allows  to conclude that $t\mapsto Y_t \diamond W^H_t$ is Bochner integrable of index $q(p_0)$.
\myqed

\smallskip

We end this section with the following lemma, the proof of which is obvious in view of Proposition \ref{idazuhdiz}, that will be useful in the proof of Theorem \ref{deftheo} below.

\begin{lemma}
\label{psdfojzoeijdozeijdoezdzoeijdzeoidjzeoijsdfpeofoierofo}
For every $p$ in $\N$, the map $(t,H) \mapsto \frac{\partial\sifbm}{\partial H}(t,H)$ is continuous from $[0,1]$ into $(({\cS}_{-p}), {\| \ \|}_{p})$. In particular, for every subset $[a,b]$ of $(0,1)$, there exists a positive real $\kappa$ such that:

\vspace{-0.25cm}

\begin{equation}
\label{sqmklqskdlmqlskd}
\forall p \in \N, \hspace{0.5cm}\underset{(s,H) \in [0,t]\times [a,b]}{\sup}{\| \tfrac{\partial\sifbm}{\partial H}(s,H)}\|_{-p} \leq \kappa.
\end{equation}

\end{lemma}

%

\section{Integral with respect to mBm through approximating fBms}
%

\hspace{0.5cm} Our aim in this section, is to construct an integral w.r.t.\h mBm using approximating integrals w.r.t. fBms. This new integral, that will be named limiting fractional Wick-It\^{o} integral, is defined at the end of this section.  The main result of this section is Theorem \ref{deftheo}, which requires the result given in Lemma \ref{ozeifjoezrifjzeroifjo} below.

Let $p_{0}$ be a fixed integer and $Y:={(Y_t)}_{t \in [0,1] }$ be an $({\cS}_{-p_0})$-valued process ({\textit{i.e.}} $Y_t$ belongs to  for every real $t$ in $[0,1]$ and $t\mapsto Y_t$ is measurable from $(0,1)$ to $({\cS}_{-p_0})$, endowed with its Borelian measure).
 As explained in Section \ref{oerifjeorifjeroifjeroifjeo}, we wish to define the integral w.r.t. an mBm $B^{h}$, noted $\int^{1}_{0} \  Y_t  \ \ d^{\diamond}B^{h}_t$, by a formula of the kind:

 \vspace{-0.25cm}

\begin{equation}
\label{gendef}
\int^{1}_{0} \  Y_t  \ \ d^{\diamond}B^{h}_t :=  \lim_{n\rightarrow\infty}  \int^1_0 \ Y_t \ \ d^{\diamond} B^{h_n}_t +  \int^{1}_{0} \ h'(t) \ Y_t \diamond \ \tfrac{\partial \sifbm}{\partial H}(t,h(t)) \ dt,
\end{equation}

where the limit holds in ${(\cS)}^*$. For this formula to make sense, it is certainly necessary that $Y$ be Bochner-integrable with respect to fBm, on $[0,1]$, of all exponents $\alpha$ in $h([0,1])$. The following technical lemma will be useful to establish Theorem \ref{deftheo} below.

\begin{lemma}
\label{ozeifjoezrifjzeroifjo}
  For any $[a,b] \subset (0,1)$ and any integer $p_0\geq 2$, there exists a positive real $\gamma_{p_0}$ such that,  for all $(t,\alpha,\alpha') \in [0,1]\times {[a,b]}^2 $,

$$ {\| W^{\alpha}_t - W^{\alpha'}_t \|}_{-p_0} \leq \gamma_{p_0}\ {|\alpha - \alpha'|}.$$
\end{lemma}

\proof
 The interval $[a,b]$ in $(0,1)$ and an integer $p_0\geq 2$ being fixed, one can write  by definition of the ${(\cS)}^*$-valued process $W$,  for all $(t,\alpha,\alpha')$ in $[0,1]\times {[a,b]}^2 $:

 \vspace{-0.3cm}

$$ {\|   W^{\alpha}_t - W^{\alpha'}_t  \|}^2_{-p_0} =  \sum^{+\infty}_{k=0}  \tfrac{ {(M_{\alpha}(e_k)(t) -M_{\alpha'}(e_k)(t) )}^2}{{(2k+2)}^{2p_0}}.$$

 Besdides, for all $(t,k)$ in $[0,1]\times\N$, the function $\alpha \mapsto M_{\alpha}(e_k)(t)$ is differentiable on $(0,1)$ (this is shown in \cite[Lemma 5.5]{JLJLV1}).
Using point 1 of \cite[Lemma 5.6]{JLJLV1} and the mean value theorem, one obtains the following fact:  there exists a positive real $\rho$ such that for all $(t,\alpha,\alpha',k) \in [0,1]\times {[a,b]}^2 \times \N$:
$$| M_{\alpha}(e_k)(t) -M_{\alpha'}(e_k)(t) |  \leq \rho \  {(k+1)}^{2/3} \ \ln(k+1) \ |\alpha - \alpha'|.$$

As a consequence, we get

 \vspace{-0.5cm}

 $${\| W^{\alpha}_t - W^{\alpha'}_t \|}^2_{-p_0} \leq {|\alpha - \alpha'|}^2 \  \rho^2 \ {\textstyle  \sum^{+\infty}_{k=0} }  \tfrac{{(k+1)}^{4/3} \ \ln^{2}(k+1)}{2^{2p_0} {(k+1)}^{2p_0}} =: |\alpha - \alpha'|^2 \  \gamma^2_{p_0}.$$

Since $p_0\geq2$,  $\gamma_{p_0}$ is finite and the proof is complete.
\myqed

\medskip

The following theorem, which constitutes the main result of this section, ensures us that  Bochner-integrability of $Y$ on $[0,1]$ is sufficient to guarantee that both the sequence ${(\int^1_0 \ Y_t \ \ d^{\diamond} B^{h_n}_t)}_{n\in \N}$   and  the quantity $\int^{1}_{0} \ h'(t) \ Y_t \diamond \ \tfrac{\partial \sifbm}{\partial H}(t,h( t)) \ dt$ exist and belong to ${(\cS)}^*$. It also establishes that the sequence ${(\int^1_0 \ Y_t \ \ d^{\diamond} B^{h_n}_t)}_{n\in \N}$ converges in ${(\cS)}^*$. For any integer $p_0$, $q(p_0)$ still denotes an integer defined as before Lemma \ref{szfojiojfdsfoifoj}

\begin{theorem}
\label{deftheo}
For any  process $Y:={(Y_t)}_{t \in [0,1] }$ that is Bochner integrable on $[0,1]$ of index $p_0$, the sequence $
{(\int^1_0 \ Y_t \ \ d^{\diamond} B^{h_n}_t)}_{n\in \N}$  and  the quantity $\int^{1}_{0} \ h'(t) \ Y_t \ \diamond \ \tfrac{\partial
\sifbm}{\partial H}(t,h(t)) \ dt$ are well defined in  ${(\cS)}^*$ and  both belong to $({\cS}_{-q(p_0)})$. Moreover the sequence ${(\int^1_0 \ Y_t \ \ d^{\diamond} B^{h_n}_t)}_{n\in \N}$  converge in $({\cS}_{-q(p_0)})$ to  $\int^1_0 \ Y_t \ \diamond \ W^{h(t)}_{t} \ dt$.
\end{theorem}

\proof
The existence of the sequence  ${(\int^1_0 \ Y_t \ \ d^{\diamond} B^{h_n}_t)}_{n\in \N}$ in  ${({\cS}_{-q(p_0)})}^{\N}$ is a straightforward consequence of Lemma \ref{szfojiojfdsfoifoj} since it has been proven there that $\int^1_0 \ Y_t \ \ d^{\diamond} B^{H}_t$ is well defined, for any $H$ in $(0,1)$, and belongs to $({\cS}_{-q(p_0)})$; $q(p_0)$ being independent from $H$. The scheme of the proof of the existence of $\int^1_0 \ Y_t \  \diamond W^{h(t)}_t \ dt$ is the same that the one we used, in the proof of Lemma  \ref{szfojiojfdsfoifoj}, to show the existence of $\int^1_0 \ Y_t \  \diamond W^{H}_t \ dt$.
One only needs to show that $\underset{t \in [0,1] }{\sup}{{\|W^{h(t)}_t\|}_{-p_0}}$ is finite for any $p_0\geq 2$. Let then $p_0\geq 2$ be fixed. Thanks to \eqref{firsozoijkqsjdefrt2}, \eqref{firsozoijkqdezdzedezddesjdefrt}  and to the upper-bound given in \cite[Theorem 3.7 point 3]{JLJLV1}, one gets:

%
%

\vspace{-0.5cm}

\begin{align*}
U_{p_0}:&=\underset{t \in [0,1] }{\sup}{{\|W^{h(t)}_t\|^2}_{-p_0}}  =  \underset{t \in [0,1] }{\sup}{ \ \sum^{+ \infty}_{k = 0}   {(M_{h(t)}({e}_k)(t))}^2  {(2k+2)}^{-2p_0} }  \\
&\leq \varrho^2_h  \sum^{+ \infty}_{k = 0} {{(2k+2)}^{-2(p_0-2/3)}},
\end{align*}

where  $\varrho_h:=\tfrac{D}{\underset{H \in h([0,1])}{\sup}{\hspace{-0.35cm}{c_H}}}$; $D$ being given in \cite[Theorem 3.7 point 3]{JLJLV1} and $c_H$ being defined right after Formula \eqref{ozeijdzoeijdzejdozeij}. Since $U_{p_0}$ is finite as soon as $p_0\geq 2$, the existence of ${\int^1_0 \ Y_t \  \diamond W^{h(t)}_t \ dt}$ is established.

 In order to show the existence of $\int^{1}_{0} \ h'(t) \ Y_t \ \diamond \ \tfrac{\partial
\sifbm}{\partial H}(t,h(t)) \ dt$, one just needs to show that  the  map $t\mapsto  h'(t)  \ Y_t\ \diamond\ \frac{\partial \sifbm}{\partial H}(t,h(t))$ is Bochner
integrable on $[0,1]$. Using the same arguments as in the proof of Lemma \ref{szfojiojfdsfoifoj} one easily prove that  $t\mapsto h'(t) \ Y_t \ \diamond \
\frac{\partial \sifbm}{\partial H}(t,h(t))$ is weakly measurable on $[0,1]$. Lemma \ref{psdfojzoeijdozeijdoezdzoeijdzeoidjzeoijsdfpeofoierofo} entails that,
\hspace{-0.3cm}$\underset{(s,H) \in \sigma_h}{\sup}{{\| \tfrac{\partial \sifbm}{\partial H}(s,H)\|}_{-p_0}} \leq \kappa$ for every $p_0$, where $\sigma_h:=[0,1]\times h([0,1])$. We hence get, 

\vspace{-0.5cm}

 $${\| h'(t) \ Y_t  \diamond  \tfrac{\partial \sifbm}{\partial H}(t,h(t))\|}_{-q(p_0)} \leq {\| Y_s \|}_{-p_0} (\underset{s \in [0,1]}{\sup}{|h'(s)|}) \underset{s \in [0,1]}{\sup}{{\|\ \hspace{-0.1cm} \tfrac{\partial \sifbm}{\partial H}(s,h(s))\|}_{-p_0}} .$$

Thus there exists $\delta \in \R^*_+$, such that $\int^1_0 {\|h'(s) \ Y_s \diamond  \frac{\partial \sifbm}{\partial H}(s,h(s))\|}_{-q(p_0)} \ ds \leq \delta \int^1_0 {\left\|Y_s \right\|}_{-p_0} \ ds <+\infty$. As a consequence,
$\int^1_0 \ h'(t) \ Y_t  \diamond  \frac{\partial \sifbm}{\partial H}(t,h(t)) \ dt$ is well defined in the sense of Bochner.

\medskip

Finally it just remains to show the convergence, in $({\cS}_{-q(p_0)})$, of ${(\int^1_0 \ Y_t \ \ d^{\diamond} B^{h_n}_t)}_{n\in \N}$ to  $\int^1_0 \ Y_t \ \diamond \ W^{h(t)}_{t} \ dt$. In view of the definition of the functions $h_n$, \eqref{eq:Xndede} and \eqref{eigfrretth2}, the equality $\int^1_0 \ Y_t \ \ d^{\diamond} B^{h_n}_t = \int^1_0 \ Y_t \ \diamond \ W^{h_n}_t \ dt$  is obvious for every $n$ in $\N$. Setting  $I_{n}:= {\|\int^1_0 \ Y_t \ \diamond \ W^{h(t)}_{t} \ dt - \int^1_0 \ Y_t \ \ d^{\diamond} B^{h_n}_t \|}_{-q(p_0)}$ and using Proposition \ref{oiefhierhiuerhieru}, \cite[Remark (2) p.92]{Kuo2} and Lemma \ref{ozeifjoezrifjzeroifjo}, one then has:

\vspace{-0.5cm}

\begin{align*}
 I_{n}&= {\| \int^1_0 \ Y_t \ \ \diamond \   (W^{h(t)}_t - W^{h_n(t)}_{t}) \ dt\|}_{\hspace{-0.1cm}-q(p_0)} \hspace{-0.25cm} \leq \int^1_0  {\|Y_t\|}_{-p_0}    {\|W^{h(t)}_t - W^{h_n(t)}_{t}\|}_{-p_0}  dt  \notag \\
  &\leq \gamma_{p_0}  \int^1_0  {\|Y_t\|}_{-p_0} \ |h(t)-h_n(t)| \ dt.
\end{align*}

The Dominated convergence theorem of Lebesgue finally allows us to write that $\lim_{n\to+\infty} I_n =0$ and thus achieves the proof.
\myqed

\begin{remark}
 The previous proof shows in particular that one does not need the pointwise convergence of ${(h_n)}_{n\in\N}$ to $h$ on the whole interval $[0,1]$ but only almost everywhere.
\end{remark}

Define the set $\Lambda_{p_0}$ by setting:

\vspace{-0.5cm}

$$\Lambda_{p_0}:=\{ {Y:=(Y_t)}_{t \in [0,1] } \in {({\cS}_{-p_0})}^{\R} : Y \text{ is Bochner integrable of index } p_0 \ \text{on } [0,1] \}.$$

\begin{corollary}
\label{szfojiojfdsfoijsidfojedzdzedzedzedzedzedzedzzsdifjsidfoj}
Let $Y$ be in $\Lambda_{p_0}$. Then the quantity
\begin{equation*}
\label{zdchzoiczoiezijdoejzeoijojijoiojiezdiojzeojioz}
 I^h_{p_0}:=  \lim_{n\rightarrow\infty} \int^1_0 \ Y_t \ d^{\diamond}B^{h_n}_t +  \int^1_0 \ h'(t)  \ Y_t \diamond \tfrac{\partial \sifbm}{\partial H}(t,h(t)) \ dt,
\end{equation*}

where the limit and the equality both hold in $({\cS}_{-q(p_0)})$, is well-defined and belongs to $({\cS}_{-q(p_0)})$. Moreover one has the equality:
\vspace{-0.25cm}

\begin{equation}
\label{zdchzoicdzedzedzedzezoiezijdoejzeoijojijoiojiezdiojzeojioz}
 I^h_{p_0} = \int^1_0 \ Y_t \diamond W^{h(t)}_t \ dt +  \int^1_0 \ h'(t)  \ Y_t \diamond \tfrac{\partial \sifbm}{\partial H}(t,h(t)) \ dt.
\end{equation}
\end{corollary}

As a consequence of Theorem \ref{deftheo} and Corollary \ref{szfojiojfdsfoijsidfojedzdzedzedzedzedzedzedzzsdifjsidfoj}
, the integral w.r.t. mBm exists as a limit of integrals w.r.t. fBms plus a second term. Thus, we are finally able to define our integral:

\begin{definition}[Limiting fractional Wick-It\^{o} integral]
\label{oidjozdee}
For any fixed integer $p_0$ and any element $Y:={(Y_t)}_{t \in [0,1] }$ of $\Lambda_{p_0}$, the integral of \hspace{0.05cm} $Y$ with respect to $B^h$ can be obtained as limits of fractional Wick-It\^{o} integral. We note $ \int^1_0 \ Y_t \ d^{\diamond}B^{h}_t$ this integral and call it  limiting fractional Wick-It\^{o} integral. It is defined by:
\vspace{-0.5cm}

\begin{equation}
\label{definale}
 \int^1_0 \ Y_t \ d^{\diamond}B^{h}_t := I^h_{p_0} =  \lim_{n\rightarrow\infty} \int^1_0 \ Y_t \ d^{\diamond}B^{h_n}_t + \int^{1}_{0} \ h'(t) \ Y_t \diamond \ \tfrac{\partial \sifbm}{\partial H}(t,h(t)) \ dt,
\end{equation}

%

\end{definition}

In view of \eqref{gendef} and of Remark \ref{pferfoeroferjfoije}, and even if $ \int^1_0 \ Y_t \ d^{\diamond}B^{h}_t$ is not only a limit of sums of fractional Wick-It\^{o} integrals,  the name limiting fractional Wick-It\^{o} integral to call $ \int^1_0 \ Y_t \ d^{\diamond}B^{h}_t$ seems to be indicate since it give the essence of it.

\begin{remark}
 The linearity of  limiting fractional Wick-It\^{o} integral as well as the equality $\int^b_a \  d^{\diamond}B^{h}_t \stackrel{{\text{a.s.}}}{=} B^{h}_b - B^{h}_a$, for any $(a,b)$ in $\R^2$ such that $a<b$, are consequences of Definition \ref{oidjozdee}.

Moreover, for any $({\cS}_{-p_0})$-valued process $Y:={(Y_t)}_{t \in [0,1] }$ that admits a  limiting fractional Wick-It\^{o} integral over a Borel subset $I$ of $\R$, if $\int_{I} X_s \ d^{\diamond}B^{h}_s$ belongs to $(L^2)$., then $\E[\int_{I} X_s \ d^{\diamond}B^{h}_s] = 0$.
\end{remark}

We shall compare the integral w.r.t. mBm, obtained in Definition  \ref{oidjozdee}, to the one defined
with the direct approach of \cite{JLJLV1}. This the goal of the next section.


\section{A comparison between multifractional Wick-It\^{o} integral and limiting fractional Wick-It\^{o} integral}
\label{comparaison}
A multifractional Wick-It\^{o} integral with respect to mBm was defined in \cite{JLJLV1}. In addition It\^{o} fromulas (in both weak and strong senses) as well as a Tanaka formula were provided. It is interesting to check whether it coincides with the one provided by Definition \ref{oidjozdee}.\h In that view, we need to adapt the definition of multifractional Wick-It\^{o} integral with respect to mBm given in \cite{JLJLV1}, which used Pettis integrals, to deal with Bochner integrals.

\subsection{Multifractional Wick-It\^{o} integral in Bochner sense}

\begin{definition}[Multifractional Wick-It\^{o} integral in Bochner sense]
\label{eigfrretoioth}
Let $I$ be a Borelian connected subset of $[0,1]$, $B^{h}:={(B^h_t)}_{t\in I}$ be a multifractional Brownian motion and $Y:={(Y_t)}_{t\in I}$ be a ${(\cS)}^*$-valued process such that:
\bit
\iti There exists $p \in \N$ such that $Y_t \in ({\cS}_{-p})$ for almost every $t \in I$,
\itii the process $t \mapsto Y_t\diamond W^{h}_t$ is Bochner integrable on $I$.
\eit

$Y$ is then said to be integrable on $I$ with respect to mBm in the Bochner sense or to admit a multifractional Wick-It\^{o} integral. This integral, noted
$\int_{I} Y_s \ dB^{h}_s \ ds$, is defined by $\int_{I} Y_s \ dB^{h}_s \ ds := \int_{I} Y_s \diamond W^{h}_s \ ds$. \end{definition}

\begin{remark}
\label{ejosijosjosjfdsoijsfddoijfd}
From the definition of  ${(W^h_t)}_{t\in[0,1]}$ \cite[Proposition $5.9$]{JLJLV1}, and the proof of Lemma \ref{szfojiojfdsfoifoj}, it is clear that every ${(\cS)}^*$-valued process $Y:={(Y_t)}_{t\in I}$ which is Bochner integrable on $I$, of index $p_0$, is integrable on $I$ with respect to mBm, in the Bochner sense. Moreover $\int_{[0,1]} Y_t\ dB^h_t$ belongs to $({\cS}_{-q(p_0)})$, where $q(p_0)$ was defined just before Lemma \ref{szfojiojfdsfoifoj}.
\end{remark}

\subsection{A comparison between multifractional Wick-It\^{o} integral and limiting fractional Wick-It\^{o} integral}

In order to compare our two integrals with respect to mBm when they both exist, it seems natural
 to assume that $Y={(Y_t)}_{t\in[0,1]}$ is a Bochner integrable process of index $p_0 \in \N$. We keep notations of the previous sections, in particular for $p_0$ and $q(p_0)$.


\begin{theorem}
\label{pfodsfpkofpssokdfpskodfp2}
Let $Y={(Y_t)}_{t\in[0,1]}$ be a Bochner integrable process of index $p_0 \in \N$. Then $Y$ is integrable with respect to mBm in both senses of Definition \ref{oidjozdee} and Definition \ref{eigfrretoioth}. Moreover  $\int_{I}  Y_t\ dB^{h}_t \ dt$ and  $\int_{[0,1]} Y_t\ d^{\diamond}B^{h}_t$ are equal in $(\cS^*)$.
\end{theorem}

{\bfseries Proof:}
The existence of both integrals  $\int_{[0,1]} Y_t\ d^{\diamond}B^{h}_t$ and $\int_{I}  Y_t\ dB^{h}_t \ dt$ is obvious in view of Theorem \ref{deftheo} and
Remark \ref{ejosijosjosjfdsoijsfddoijfd}.\h Moreover Equalities \eqref{definale}, \eqref{zdchzoicdzedzedzedzezoiezijdoejzeoijojijoiojiezdiojzeojioz} and \eqref{firsozoijkqdezdzedezddesjdefrt} allow us to write:

\vspace{-0.5cm}

\begin{align*}
\int^1_0 \ Y_t \ d^{\diamond}B^{h}_t &= \int^1_0 \ \left( Y_t \diamond W^{h(t)}_t +  h'(t)  \ Y_t \diamond \tfrac{\partial \sifbm}{\partial H}(t,h(t)) \right) \ dt = \int^1_0 \ Y_t \diamond W^{h}_t  \ dt.\\
&= \int^1_0 \ Y_t \ dB^{h}_t. 
\end{align*}

\vspace{-0.9cm}
\myqed

 \section*{Acknowledgments} The author thanks Jacques L\'{e}vy V\'{e}hel for many helpful remarks and comments and the anoymous referee for his suggestions that highly improved Sections 4 and 5 of this paper.



 { \footnotesize
\medskip
\medskip
 \vspace*{1mm}

\noindent {\it Joachim Lebovits}\\
Mathematic Center of Heidelberg \& University of Heidelberg\\
Im Neuenheimer Feld, 294, 69120 Heidelberg, Germany\\
E-mail: {\tt \url{joachim.lebovits@uni-heidelberg.de}}\\ \\

 \end{document}